\documentclass[12pt]{amsart}
\usepackage{amsthm, amssymb, amsmath, amsxtra, latexsym}
\usepackage[all]{xy}
\usepackage{enumerate}

\theoremstyle{plain}
\newtheorem{theorem}{Theorem}[section]
\newtheorem{lemma}[theorem]{Lemma}
\newtheorem{corollary}[theorem]{Corollary}
\newtheorem{proposition}[theorem]{Proposition}
\theoremstyle{remark}
\newtheorem{remark}[theorem]{Remark}

\newtheorem{example}[theorem]{Example}

\theoremstyle{definition}
\newtheorem{definition}[theorem]{Definition}
\theoremstyle{remark}

\begin{document}

\begin{title}
[Constructions of $E_{\mathcal{VC}}$ and $E_{\mathcal{FBC}}$ for  CAT(0) groups]
{Constructions of $E_{\mathcal{VC}}$ and $E_{\mathcal{FBC}}$ for groups acting on CAT(0) spaces} 
\author[D.Farley]{Daniel Farley}
      \address{Department of Mathematics and Statistics\\
               Miami University \\
               Oxford, OH 45056\\
               }
      \email{farleyds@muohio.edu}

\end{title}

\begin{abstract}
If $\Gamma$ is a group acting properly by semisimple isometries on a proper CAT(0) space $X$, then we build models for the classifying spaces 
$E_{\mathcal{VC}}\Gamma$ and $E_{\mathcal{FBC}}\Gamma$ under the additional assumption that the action of $\Gamma$ 
has a well-behaved collection of axes in $X$.  We conjecture that the latter hypothesis is satisfied in a large range of cases.

Our classifying spaces are natural variations of the constructions due to Connolly, Fehrman, and Hartglass \cite{CFH} of $E_{\mathcal{VC}}\Gamma$ for crystallographic
groups $\Gamma$.          
\end{abstract}

                                                                                

\maketitle

\section{Introduction}

We say that a collection $\mathcal{F}$ of subgroups of $\Gamma$ is a \emph{family} if $\mathcal{F}$ is closed under
conjugation and passage to subgroups.  If $\Gamma$ is a discrete group, then $E_{\mathcal{F}}\Gamma$, the \emph{classifying complex of $\Gamma$ with isotropy in
$\mathcal{F}$}, is a $\Gamma$-CW complex such that:
\begin{enumerate}
\item if $c$ is a cell of $E_{\mathcal{F}}\Gamma$ and $\gamma \in \Gamma$ leaves $c$ invariant as a set,  then $\gamma$ fixes $c$ pointwise;
\item if $H \in \mathcal{F}$, then the set of points fixed by $H$ is contractible;
\item if $H \leq \Gamma$ and $H \not \in \mathcal{F}$, then the set of points fixed by $H$ is empty.
\end{enumerate}
Let $\mathcal{VC}$ denote the family of virtually cyclic subgroups of $\Gamma$.  
The classifying space $E_{\mathcal{VC}}\Gamma$ can be used to help compute the algebraic $K$-theory of the group $\Gamma$ if the Farrell-Jones isomorphism
conjecture has been proved for $\Gamma$.  Recent work by several authors (see, for instance, \cite{Khan}) 
shows that one can use the space $E_{\mathcal{FBC}}\Gamma$ for similar purposes.
Here $\mathcal{FBC}$ is the family of finite-by-cyclic subgroups of $\Gamma$, namely those that map onto a cyclic group with finite kernel. 

Connolly, Fehrman, and Hartglass \cite{CFH} built classifying spaces $E_{\mathcal{VC}}\Gamma$ for crystallographic groups $\Gamma$.  Their construction can be summarized as follows.
We let $\Gamma$ be an $n$-dimensional crystallographic group. Thus, $\Gamma$ acts properly and cocompactly by isometries on Euclidean $n$-space $\mathbb{R}^{n}$.  
Each element $\gamma \in \Gamma$ will either fix a point in $\mathbb{R}^{n}$ or act by translation on some line $\ell \subseteq \mathbb{R}^{n}$.  In the latter case, 
we call such a line $\ell$ an \emph{axis} for $\gamma$.  
Given an axis $\ell$, we let $\mathbb{R}^{n}_{\ell}$ be the set of all lines 
in $\mathbb{R}^{n}$ that are parallel to $\ell$.  It is rather clear that $\mathbb{R}^{n}_{\ell}$ is 
naturally isometric to $\mathbb{R}^{n-1}$.  We let $f_{\ell}: \mathbb{R}^{n} \rightarrow \mathbb{R}^{n-1}_{\ell}$
be the quotient map, and let $M(f_{\ell})$ be the mapping cylinder of $f_{\ell}$.  Now we let $\ell$ range over all possible axes of elements $\gamma \in \Gamma$,
and we glue together all of the mapping cylinders $M(f_{\ell})$ together along the tops, which are identical copies of $\mathbb{R}^{n}$.  
The resulting space (with the cell structure described in \cite{CFH} and with respect to a natural $\Gamma$-action) is $E_{\mathcal{VC}}\Gamma$.

The Connolly-Fehrman-Hartglass construction is particularly useful in calculations.  Since the Farrell-Jones isomorphism conjecture is known for crystallographic
groups \cite{FJ}, one can use a theorem of Bartels \cite{Bartels} to conclude that the algebraic $K$-theory of
$\Gamma$ is a direct sum of two pieces:  one contributed by the finite subgroups, and the other ($Nil$) part contributed by the
infinite virtually cyclic subgroups of $\Gamma$.  The contribution of finite subgroups can be computed from $E_{\mathcal{FIN}}\Gamma$, where $\mathcal{FIN}$ is
the family of finite subgroups.  The contribution of infinite virtually cyclic subgroups comes from the subcomplex of $E_{\mathcal{VC}}\Gamma$ having infinite virtually
cyclic cell stabilizers, and such cell stabilizers occur only at the bases of the mapping cylinders.  In other words (to summarize), the algebraic $K$-theory 
is the direct sum of two pieces: one is computed from the common top of the mapping cylinders (which is homeomorphic to $\mathbb{R}^{n}$), 
and the other is computed from the bases of the mapping cylinders (each of which is homeomorphic to $\mathbb{R}^{n-1}$).                                
       
The goal of this paper is to show that a version of the construction from \cite{CFH} is available for a wide variety of groups $\Gamma$ acting by isometries on CAT(0) spaces $X$.  Let $X$ be a proper CAT(0) space.  Suppose that $\Gamma$ acts properly on $X$ by \emph{semisimple} isometries:  that is, every $\gamma \in \Gamma$ either
fixes a point in $X$ or acts on a line $\ell$ (again called an \emph{axis}) by translation.  For instance, if $\Gamma$ acts properly and cocompactly by isometries
on $X$, then the action of $\Gamma$ is automatically by semisimple isometries (\cite{BH}, Proposition 6.10(2), page 233).  
We let $\mathcal{A}(X; \Gamma)$ denote the \emph{space of axes}; thus, $\mathcal{A}(X; \Gamma)$ is the set of all lines $\ell \subseteq X$ for which there is $\gamma \in \Gamma$ acting by translations on $\ell$.  We write $d(\ell_{1}, \ell_{2}) = K$
if the axes $\ell_{1}$ and $\ell_{2}$ are parallel and bound a flat strip of width $K$; we write $d(\ell_{1}, \ell_{2}) = \infty$ if $\ell_{1}$ and $\ell_{2}$ are not
parallel.  The function $d$ is a metric on $\mathcal{A}(X; \Gamma)$, and makes $\mathcal{A}(X; \Gamma)$ into a disjoint union of CAT(0) spaces.  The space $E_{\mathcal{VC}}\Gamma$ that we build in this paper is (roughly) the join $X \ast \mathcal{A}(X; \Gamma)$, although there are complications as we now explain.  (Note
that we use a join construction rather than a mapping cylinder construction.  The join is easier to work with, and it doesn't make $K$-theory
computations more complicated.  The latter observation essentially follows from the theorem of Bartels \cite{Bartels}.)  

The main difficulty comes from the requirement that $E_{\mathcal{VC}}\Gamma$ have a CW structure.  There is no reason, in general, to expect that the CAT(0) space $X$
(let alone $\mathcal{A}(X; \Gamma)$)
will have such a structure, so we are forced to produce our own.  Our method is to replace both $X$ and $\mathcal{A}(X; \Gamma)$ by nerves of suitable covers.  Our
procedure for producing these covers requires that $\mathcal{A}(X; \Gamma)$ be a disjoint union of proper CAT(0) spaces.  Unfortunately, this is not always the case,
even if the action of $\Gamma$ is cocompact, so we add a hypothesis.  

We need some definitions first for background.  
We let $\mathcal{L}(X)$ denote the space of all lines in $X$.  The space $\mathcal{L}(X)$
is given the metric $d$ (as defined above for $\mathcal{A}(X; \Gamma)$).  This metric makes $\mathcal{L}(X)$ a disjoint union of proper CAT(0) spaces.   We say that a subspace $K$ of $\mathcal{L}(X)$ is \emph{component-wise convex} if the intersection of $K$ with each connected component of $\mathcal{L}(X)$
is convex in the ordinary sense.  

We can now describe the hypothesis that we need.  If $\mathcal{A} \subseteq \mathcal{A}(X; \Gamma)$, we say that $\mathcal{A}$ is \emph{a $\Gamma$-well-behaved space of axes} if $\mathcal{A}$ is closed and component-wise convex as a subset of $\mathcal{L}(X)$, $\Gamma \cdot \mathcal{A} = \mathcal{A}$, and, for any $\ell \in \mathcal{A}(X; \Gamma)$, 
$\mathcal{A}$ contains at least one axis parallel to $\ell$.  (In fact, the space $\mathcal{A}(X; \Gamma)$ satisfies all of these conditions, except that it need not be closed.)  We have the following theorem: 
\begin{theorem}
If $\Gamma$ acts properly by semisimple isometries on a proper CAT(0) space $X$ and there is some $\Gamma$-well-behaved space of
axes $A$, then there are Connolly-Fehrman-Hartglass constructions of $E_{\mathcal{VC}}\Gamma$ and $E_{\mathcal{FBC}}\Gamma$ based on the action of $\Gamma$ on $X$.
\end{theorem} 
This statement of our main theorem is a place holder.  The precise statements appear in Theorem \ref{main} and Corollary \ref{FBC}.  
The $E_{\mathcal{VC}}\Gamma$ we build is (roughly) $X \ast \mathcal{A}$,
although in practice we must replace both spaces with suitable nerves of covers.  
One builds the $E_{\mathcal{FBC}}\Gamma$ complexes by using a slight modification of our techniques for building $E_{\mathcal{VC}}\Gamma$ complexes. 
           
The hypothesis that there exists a $\Gamma$-well-behaved space of axes is clearly undesirable, but I believe that the simplicity of the construction compensates for this
flaw somewhat.  Moreover, it seems reasonable to guess that the hypothesis will be satisfied in the great majority of cases.  In Example \ref{silly} we give an 
example of a space $X$ and a group $\Gamma$ such that $\mathcal{A}(X; \Gamma)$ is not a $\Gamma$-well-behaved space
of axes.  The example is rather artificial, and it is unclear whether such badly behaved examples occur ``naturally".                
The construction of the main theorem also
appears to be economical and useful for computations (like the construction of $E_{\mathcal{VC}}\Gamma$ for crystallographic groups from \cite{CFH}). 

While this paper was being completed, Wolfgang Lueck \cite{Lueck}
described a construction of $E_{\mathcal{VC}}\Gamma$ for any group $\Gamma$ acting properly and cocompactly
on a CAT(0) space $X$.  He is able to build classifying spaces for all CAT(0) groups (without our additional hypothesis) and he gives
bounds for the topological dimension of his construction.  His construction is also different from ours.  I hope that readers will see
virtues in both approaches.
  
I would like to thank Qayum Khan for suggesting the problem of building complexes $E_{\mathcal{FBC}}\Gamma$ to me while I was visiting Vanderbilt in 2008.  I
thank Ross Geoghegan for suggesting the method of using nerves of covers.

\section{A General Construction}

\begin{definition} \label{family}
If $\Gamma$ is any group, then a \emph{family} $\mathcal{F}$ of subgroups
of $\Gamma$ is a collection of subgroups that is closed under conjugation by elements of $\Gamma$
and passage to subgroups.  If $\Gamma$ is any group, then we let $\mathcal{FIN}$, $\mathcal{VC}$, and $\mathcal{FBC}$
denote (respectively) the families of finite, virtually cyclic, and \emph{finite-by-cyclic groups}, i.e., groups which map onto a cyclic group with finite
kernel.  We  will also let $\mathcal{VC}_{\infty}$ denote the difference $\mathcal{VC} - \mathcal{FIN}$ and $\mathcal{FBC}_{\infty}$ denote $\mathcal{FBC}-\mathcal{FIN}$.
(Note that neither $\mathcal{VC}_{\infty}$ nor $\mathcal{FBC}_{\infty}$ is a family, since neither collection contains the trivial group.)
\end{definition}

\begin{definition} \label{cx}
Let $X$ be a $\Gamma$-CW complex.  Suppose that, if $c \subseteq X$ is a cell of $X$, then $\gamma \cdot c = c$ if and only if
$\gamma$ fixes $c$ pointwise.  Let $\mathcal{F}$ be a family of subgroups of $\Gamma$.  We say that $X$ is an \emph{$E_{\mathcal{F}}\Gamma$-complex} if
\begin{enumerate}
\item $X$ is contractible;
\item whenever $H \in \mathcal{F}$, the fixed set $Fix(H) = \{ x \in X \mid \gamma \cdot x = x ~for~all~ \gamma \in H \}$ is contractible;
\item whenever $H \not \in \mathcal{F}$, $Fix(H)$ is empty.
\end{enumerate}
Let $\mathcal{F}'$ be a family of subgroups of $\Gamma$ containing $\mathcal{F}$.  We say that $X$ is an \emph{$I_{\mathcal{F}' - \mathcal{F}}\Gamma$-complex} if
\begin{enumerate}
\item whenever $H \in \mathcal{F}' - \mathcal{F}$, $Fix(H)$ is contractible;
\item whenever $H \not \in \mathcal{F}'$, $Fix(H)$ is empty.
\end{enumerate}  
\end{definition}
 
\begin{remark}  
Note that the complex $I_{\mathcal{F}'-\mathcal{F}}\Gamma$ isn't necessarily contractible.  In our applications, $I_{\mathcal{F}'-\mathcal{F}}\Gamma$
won't even be connected.
\end{remark} 
 
\begin{proposition} \label{promoting}
Let $\Gamma$ be a group.  Let $\mathcal{F} \subseteq \mathcal{F}'$ be
families of subgroups of $\Gamma$.  The join
$$ \left( E_{\mathcal{F}}\Gamma \right) \ast \left( I_{\mathcal{F}' - \mathcal{F}}\Gamma \right)$$
is an $E_{\mathcal{F}'}\Gamma$-complex.
\end{proposition}

\begin{proof}
For any two topological spaces $X$ and $Y$ with $\Gamma$-actions, we make the following two observations:  
(i) the join $X \ast Y$ is contractible if at least one of the spaces $X$ and $Y$
is contractible; (ii) $Fix_{X \ast Y}(H) = Fix_{X}(H) \ast Fix_{Y}(H)$, for any $H \leq \Gamma$.

We let $Z = (E_{\mathcal{F}}\Gamma) \ast (I_{\mathcal{F}' - \mathcal{F}}\Gamma)$ and give $Z$ the natural cell structure. 
We need to show that $Z$ satisfies (1)-(3) with respect to the family $\mathcal{F}'$.  The contractibility of $Z$
follows from (i), since $E_{\mathcal{F}}\Gamma$ is contractible.  Let $H \leq \Gamma$ satisfy $H \not \in \mathcal{F}'$.
Since $Fix_{E_{\mathcal{F}}\Gamma}(H) = \emptyset = Fix_{I_{\mathcal{F}' - \mathcal{F}}\Gamma}(H)$, it follows from (ii) that
$Fix_{Z}(H) = \emptyset$.

We now need to check (2).  Suppose first that $H \in \mathcal{F}' - \mathcal{F}$.  Since $Fix_{E_{\mathcal{F}}\Gamma}(H) = \emptyset$ and
$Fix_{I_{\mathcal{F}' - \mathcal{F}}\Gamma}(H)$ is contractible, $Fix_{Z}(H)$ is contractible by (ii).  Suppose finally that $H \in \mathcal{F}$.
Since $Fix_{E_{\mathcal{F}}\Gamma}(H)$ is contractible, $Fix_{Z}(H)$ is contractible by (ii).

Finally, we note that, for any cell $c \subseteq Z$, an element $\gamma \in \Gamma$ leaves $c$ invariant if and only if $\gamma$ fixes $c$ pointwise, 
since the cells of $E_{\mathcal{F}}\Gamma$ and $I_{\mathcal{F}' - \mathcal{F}}\Gamma$ have this property.          
\end{proof}

\section{Preliminaries about CAT($0$) Spaces}

All of the material in this section is either taken directly from \cite{BH} or closely based on that source.

\subsection{General Definitions}

\begin{definition} \label{basic}
Let $X$ be a metric space.  We say that $X$ is a \emph{geodesic metric space} if, for any $x$, $y \in X$, there
is a map, called a \emph{geodesic segment}, $\ell_{x,y}: [0, d(x,y)] \rightarrow X$ such that, 
$\ell_{x,y}(0) = x$, $\ell_{x,y}(d(x,y)) = y$, and for any $t_{1}$, $t_{2} \in [0, d(x,y)]$,
$d(\ell_{x,y}(t_{1}), \ell_{x,y}(t_{2})) = |t_{1} - t_{2}|$.  We will also call the image of such a map
a geodesic segment.  Similarly, a \emph{geodesic line} is an isometric embedding $\ell : \mathbb{R} \rightarrow X$.
We will also call the image of such an embedding a geodesic line.  

Now suppose that $X$ is a geodesic metric space.  Let $x,y, z \in X$.  Choose geodesic segments $[x,y]$, $[y,z]$,
and $[x,z]$ (these need not be unique, a priori).  We choose three points $\bar{x}, \bar{y}, \bar{z} \in \mathbb{E}^{2}$
such that $d_{X}(x,y) = d_{\mathbb{E}^{2}}(\bar{x},\bar{y})$, $d_{X}(y,z) = d_{\mathbb{E}^{2}}(\bar{y},\bar{z})$, and
$d_{X}(x,z) = d_{\mathbb{E}^{2}}(\bar{x},\bar{z})$.  The segments $[x,y]$, $[y,z]$, and $[x,z]$ are in natural isometric one-to-one
correspondence with the segments $[\bar{x},\bar{y}]$, $[\bar{y},\bar{z}]$, and $[\bar{x},\bar{z}]$ (respectively).  If 
$a \in [x,y]$, $[y,z]$, or $[x,z]$, then the corresponding point $\bar{a}$ on (respectively) 
$[\bar{x},\bar{y}]$, $[\bar{y},\bar{z}]$, or $[\bar{x},\bar{z}]$ is called a \emph{comparison point}.  We say that the
geodesic triangle $\Delta = [x,y] \cup [y,z] \cup [x,z]$ satisfies the \emph{CAT(0) condition} if, for any pair of points $a$, $b \in
\Delta$ and pair of comparison points $\bar{a}$, $\bar{b}$, we have
$$ d_{X}(a,b) \leq d_{\mathbb{E}^{2}}(\bar{a}, \bar{b}).$$
The geodesic metric space $X$ is \emph{CAT(0)} if every geodesic triangle satisfies the CAT(0) condition.    

We say that a metric space $X$ is \emph{proper} if each closed metric ball 
$\overline{B}_{\epsilon}(x) = \{ y \in X \mid d(x,y) \leq \epsilon \}$ ($0 < \epsilon < \infty$) is compact. 
\end{definition}

\begin{remark}  
It will be helpful if we allow two points in a metric space to be an infinite distance apart.  We will always allow
this in what follows.  In particular, if $X$ is a disjoint union of metric spaces $X_{i}$, then we metrize $X$ by letting $d_{X}$ agree
with the metrics on the $X_{i}$, and setting $d_{X}(x,y) = \infty$ if $x \in X_{i}$ and $y \in X_{j}$ for $i \neq j$.
\end{remark}

\begin{definition}  
If $X$ is a geodesic metric space (respectively, if $X$ is a disjoint union of geodesic metric spaces), then we say that
a subset $K \subseteq X$ is \emph{convex} (respectively, \emph{component-wise convex})
 if, for any two points $x$, $y \in K$ satisfying $d(x,y) < \infty$, every geodesic
segment connecting $x$ to $y$ is contained in $K$.
\end{definition}

\begin{definition}
If $\ell_{1}: \mathbb{R} \rightarrow X$ and $\ell_{2}: \mathbb{R} \rightarrow X$ are geodesic lines in a complete CAT($0$) space $X$, then
we say that $\ell_{1}$ and $\ell_{2}$ are \emph{parallel} (or \emph{asymptotic}) if there is some constant $K$ such that 
$d_{X}(\ell_{1}(t), \ell_{2}(t)) \leq K$, for all $t \in \mathbb{R}$.
\end{definition} 

\begin{proposition} \label{dump}
Let $X$ be a complete CAT($0$) space.
\begin{enumerate}
\item Any two points in $X$ can be connected by a unique geodesic (i.e., $X$ is \emph{uniquely geodesic}).
\item $X$ is contractible.  Every open or closed ball of finite radius in $X$ is convex, and therefore CAT($0$) and contractible.
\item Every bounded subset $B$ in $X$ has a unique center $c$.  If an isometry $\gamma$ leaves the set $B$ invariant, then 
$\gamma$ fixes $c$.  If $B$ is a ball and $\gamma$ fixes $c$, then $B$ is invariant under the action of $\gamma$.
\item (Flat Strip Theorem)  If $\ell_{1}$ and $\ell_{2}$ are parallel geodesic lines in $X$, then, for some $K \geq 0$, there is
an isometric embedding $\rho: [0,K] \times \mathbb{R} \rightarrow X$, where $\ell_{1} = \rho( \{ 0 \} \times \mathbb{R})$ and
$\ell_{2} = \rho( \{ K \} \times \mathbb{R})$. 
\end{enumerate}  
\begin{proof}
These facts are proved on pages 160, 161, 179, and 182 of \cite{BH} (respectively).
\end{proof}

\end{proposition}

\subsection{Spaces of Lines and Axes}

\begin{definition} \label{sol}
Let $X$ be a proper CAT(0) space endowed with a proper $\Gamma$-action by semisimple isometries.
The \emph{space of lines} in $X$, denoted $\mathcal{L}(X)$, is defined as follows:
$$ \mathcal{L}(X) = \{ \ell \subseteq X \mid \ell ~is~a~geodesic~line \}.$$
If $\ell_{1}, \ell_{2} \in \mathcal{L}(X)$, then we write
$$ d(\ell_{1}, \ell_{2}) = \left\{ \begin{array}{ll} k & ~if~ \ell_{1} ~and~ \ell_{2} ~bound~a~flat~strip~of~width~ k \\
\infty & ~if~ \ell_{1} ~and~ \ell_{2} ~are~not~parallel. \end{array} \right.$$
The given assignment $d$ is a function by the Flat Strip Theorem (Proposition \ref{dump}(4)).  
It follows easily from the Product Decomposition Theorem (Theorem \ref{decomp}, below) that $d$ is a metric.

We say that a line $\ell \subseteq X$ is an \emph{axis} for $\gamma \in \Gamma$ if $\gamma$ acts on $\ell$ by translation.
The \emph{space of axes} for elements of $\Gamma$ in $X$, denoted $\mathcal{A}(X; \Gamma)$, is defined as follows:
$$\mathcal{A}(X; \Gamma) = \{ \ell \in \mathcal{L}(X) : \ell ~is~an~axis~for~some~ \gamma \in \Gamma \}.$$
\end{definition}

\begin{theorem} \label{decomp}
(Product Decomposition Theorem) Let $X$ be a proper CAT(0) space and let $c: \mathbb{R} \rightarrow X$ be a geodesic line.
\begin{enumerate}
\item The union of the images of all geodesic lines $c': \mathbb{R} \rightarrow X$ parallel to $c$ is a closed convex subspace
$X_{c}$ of $X$.
\item Let $p$ be the restriction to $X_{c}$ of the projection from $X$ to the complete convex subspace $c(\mathbb{R})$.
Let $X^{0}_{c} = p^{-1}(c(0))$.  Then, $X_{c}^{0}$ is closed and convex (in particular, it is a proper CAT(0) space) and
$X_{c}$ is canonically isometric to the product $X_{c}^{0} \times \mathbb{R}$.
\item Any (image of a) geodesic line parallel to $c$ has the form $\{ x \} \times \mathbb{R}$ for some $x$ in  $X_{c}^{0}$.
(Conversely, any subspace of the form $\{ x \} \times \mathbb{R}$ is the image of a geodesic line parallel to $c$.)
\end{enumerate}
\end{theorem}

\begin{proof}
This  theorem is exactly the same as the Product Decomposition Theorem from \cite{BH} (page 183), except that the latter theorem contains neither the word ``proper",
nor the word ``closed", nor statement (3).  Statement (3) follows directly from the correspondence $j$ as defined on page 183 of \cite{BH}.  
(In fact, (3) is the canonical 
identification which was alluded to in (2).)  Therefore, it is enough
to prove (1) and (2).     

Statement (1) is proved in \cite{BH} on page 183, except for the statement that $X_{c}$ is closed.  The latter 
follows easily from the fact that $X$ is proper and from Ascoli's Theorem
from page 290 of \cite{Munk}.  The details of the argument are routine, and are omitted.  Statement (2) follows directly from (1) and the Product Decomposition Theorem
from \cite{BH}.

\end{proof}

\begin{proposition} \label{lines}
If $X$ is a proper CAT(0) space and $\Gamma$ acts by isometries on $X$, 
then the space of lines $\mathcal{L}(X)$ is a disjoint union of proper CAT(0) spaces on which $\Gamma$ acts isometrically. 
\end{proposition}

\begin{proof}
It is clear that $\Gamma$ acts isometrically on $\mathcal{L}(X)$.  We need to show that $\mathcal{L}(X)$ is a disjoint union of
proper CAT(0) spaces.

Fix a geodesic line $c: \mathbb{R} \rightarrow X$.  We consider $X_{c}$, the union of all lines parallel to $c$.  By Theorem \ref{decomp}(2)
$X_{c} = X_{c}^{0} \times \mathbb{R}$.  Let $\mathcal{L}_{c}(X) = \{ \ell \in \mathcal{L}(X) \mid \ell ~is~parallel~to~ c(\mathbb{R}) \}$.
The space $\mathcal{L}_{c}(X)$ is naturally isometric to $X_{c}^{0}$ by Theorem \ref{decomp}(3).  Therefore, $\mathcal{L}_{c}(X)$ is proper by
Theorem \ref{decomp}(2).  The space $\mathcal{L}(X)$ is the disjoint union
of all $\mathcal{L}_{c}(X)$, as $c$ ranges over a maximal collection of pairwise non-parallel lines.
\end{proof}

The following Proposition isn't needed in order to understand the main construction of this paper.  It does show, however, that the space of axes
comes very close to being ``$\Gamma$-well-behaved" (see Definition \ref{wb}).  Example \ref{silly} shows how the space of axes can fail to be well-behaved.
 
\begin{proposition} \label{axes}
If $\Gamma$ acts properly by semisimple isometries on the proper CAT(0) space $X$, then $\Gamma$ acts discretely by isometries on
$\mathcal{A}(X; \Gamma)$.  The space $\mathcal{A}(X; \Gamma)$ is a $\Gamma$-equivariant, component-wise convex subset of $\mathcal{L}(X)$.
\end{proposition}

\begin{proof}
Let $\ell \in \mathcal{A}(X; \Gamma)$.  By definition, there is some $\gamma \in \Gamma$ such that $\gamma$ acts on $\ell$ by translations.
Now choose an arbitrary $\gamma_{1} \in \Gamma$.  We claim that $\gamma_{1} \cdot \ell$ is an axis for $\gamma_{1} \gamma \gamma_{1}^{-1}$.
$$ (\gamma_{1}\gamma \gamma_{1}^{-1}) \cdot \gamma_{1} \cdot \ell = \gamma_{1}\gamma \ell = \gamma_{1} \ell.$$
Therefore $\gamma_{1}\ell$ is invariant under the action of $\gamma_{1}\gamma\gamma_{1}^{-1}$.  Since $\gamma_{1}\gamma\gamma_{1}^{-1}$ has infinite
order and the action of $\Gamma$ is proper, $\gamma_{1}\gamma\gamma_{1}^{-1}$ must act without a fixed point, and  therefore $\gamma_{1}\gamma\gamma_{1}^{-1}$
acts by translation on $\gamma_{1} \ell$.  This shows that $\Gamma$ acts on $\mathcal{A}(X; \Gamma)$.  The action is by isometries since $\Gamma$ 
acts on $\mathcal{L}(X)$ by isometries.

We now need to show that the action of $\Gamma$ on $\mathcal{A}(X;\Gamma)$ has discrete orbits.  Suppose otherwise.  There must exist a sequence $\ell_{1}, \ldots,
\ell_{n},  \ldots$ of axes (all in the same orbit) converging to  another axis $\ell$ (in the same orbit as the others).  We can assume that any two of the axes
are parallel, and that $d(\ell_{n}, \ell)$ is a strictly decreasing sequence, after passing to a subsequence if necessary.  It follows that $\ell,  \ell_{1},  \ldots,
\ell_{n}, \ldots$ all lie in the same component $\mathcal{L}_{c}(X)$ of $\mathcal{L}(X)$.  By Theorem \ref{decomp}(2), 
$\mathcal{L}_{c}(X)$ is isometric to $Y \times \mathbb{R}$,
where $Y$ is a proper CAT(0) space.  By Theorem \ref{decomp}(3), $\ell, \ell_{1}, \ldots, \ell_{n}, \ldots$ are identified with subspaces $\{ y \} \times \mathbb{R},$
$\{ y_{1} \} \times \mathbb{R},$ $\ldots, \{ y_{n} \} \times \mathbb{R}, \ldots$ respectively.

We consider the point $(y,0) \in \ell$.  Since $\ell$ is an axis, there is some isometry $\gamma$ that acts by translation on $\ell = \{ y \} \times \mathbb{R}$. 
Since $\gamma$ acts on $Y \times \mathbb{R}$ and preserves the line $\ell$, Proposition 5.3(4) from page 56 of \cite{BH} implies that
$\gamma_{\mid Y \times \mathbb{R}}
= \gamma_{1} \times \gamma_{2}$, where $\gamma_{1}: Y \rightarrow Y$,
$\gamma_{2}: \mathbb{R} \rightarrow \mathbb{R}$ are isometries, $\gamma_{1}$ fixes $y$, and $\gamma_{2}(t) = t + K$ for some $K \in \mathbb{R}$.

Since $\ell,$
$\ell_{1}, \ldots, \ell_{n}, \ldots$ are all in the same orbit under $\Gamma$, there is a sequence $(y_{1}, t_{1}),$  $(y_{2}, t_{2}), \ldots, (y_{n}, t_{n}),
\ldots$ of points in the orbit of $(y,0)$ under $\Gamma$, where $(y_{i}, t_{i}) \in \ell_{i}$, for all $i \in \mathbb{N}$.  Note that  
$d(\ell_{i}, \ell) = d(y_{i}, y)$, so $d_{Y}(y_{i}, y)$ is a decreasing sequence.  We can find integral powers $m_{1}, m_{2}, \ldots, m_{n}, \ldots$ where
$|\gamma_{2}^{m_{i}}(t_{i})| \leq K$, for each $i \in \mathbb{N}$.
Now consider the points
$$ (\gamma_{1}^{m_{1}}(y_{1}), \gamma_{2}^{m_{1}}(t_{1})), (\gamma_{1}^{m_{2}}(y_{2}), \gamma_{2}^{m_{2}}(t_{2})), \ldots.$$
All of these points are in the orbit $(y,0)$ according to our assumption.  Moreover, they form a bounded set, since 
$$ d((\gamma_{1}^{m_{i}}(y_{i}), \gamma_{2}^{m_{i}}(t_{i})), (y,0)) = d((\gamma_{1}^{m_{i}}(y_{i}), \gamma_{2}^{m_{i}}(t_{i})), (\gamma_{1}^{m_{i}}(y),0))
\leq \sqrt{d_{Y}(y_{i},y)^{2} + K^{2}}.$$
Finally, we note that all of the points $(\gamma_{1}^{m_{i}}(y_{1}),$ $\gamma_{2}^{m_{1}}(t_{1})), \ldots$ are distinct since the points
$\gamma_{1}^{m_{1}}(y_{1}), \gamma_{1}^{m_{2}}(y_{2}), \ldots$ are all distinct.  (Indeed,  the distances $d_{Y}(\gamma_{1}^{m_{i}}(y_{i}), y) = d_{Y}(y_{i},y)$
are all distinct.)

We've found an infinite, bounded collection of points in a single orbit.  This  violates properness of the action of $\Gamma$.  It follows that $\Gamma$ acts on
$\mathcal{A}(X;\Gamma)$ with discrete orbits, as claimed.

Lastly, we need to show that if  $\ell_{1}, \ell_{2} \in \mathcal{A}(X;\Gamma)$ and $[ \ell_{1}, \ell_{2}]$ is a geodesic segment connecting $\ell_{1}$ to $\ell_{2}$
in $\mathcal{L}(X)$, then $[\ell_{1}, \ell_{2}] \subseteq \mathcal{A}(X;\Gamma)$.  We choose a point $\ell_{3} \in [ \ell_{1}, \ell_{2}]$.  We will show first that
the orbit of $\ell_{3}$ is discrete, and then that $\ell_{3}$ is an axis.  Suppose, for a contradiction, that $\ell_{3}$ has a non-discrete orbit.  Let 
$\gamma_{1}, \gamma_{2}, \ldots, \gamma_{n}, \ldots$ be a sequence of isometries of $X$ such that $\gamma_{1} \cdot \ell_{3}, \gamma_{2} \cdot \ell_{3}, \ldots$
are all distinct and lie inside a ball of radius $1$ in $\mathcal{L}(X)$.  Let us call the center of this ball $\ell_{c}$.  The lines $\gamma_{1}\cdot \ell_{1},
\gamma_{2}\cdot \ell_{1}, \ldots$ all lie within $d(\ell_{1},\ell_{2}) + 1$ units of $\ell_{c}$.  Since the component of $\mathcal{L}(X)$ containing $\ell_{c}$ is
proper and $\Gamma$ acts discretely on $\mathcal{A}(X; \Gamma)$, the set $L_{1} = \{ \gamma_{1} \cdot \ell_{1}, \gamma_{2} \cdot \ell_{1}, \ldots, \gamma_{n} \cdot \ell_{1}, \ldots \}$ is finite.  The same reasoning applies
to the set $L_{2} = \{ \gamma_{1} \cdot \ell_{2}, \gamma_{2} \cdot \ell_{2}, \ldots, \gamma_{n} \cdot \ell_{2}, \ldots \}$, so it is finite as well.
Consider the points $\ell \in \mathcal{L}(X)$ with the following properties:
\begin{enumerate}
\item $\ell$ lies on a geodesic $[ \hat{\ell}_{1}, \hat{\ell}_{2} ]$, where $\hat{\ell}_{i} \in L_{i}$ for $i= 1,2$.
\item $d_{\mathcal{L}(X)}(\ell, \ell_{i}) = d_{\mathcal{L}(X)}(\ell_{3}, \ell_{i})$ for $i=1,2$.
\end{enumerate}
The set of all such $\ell$ is finite (indeed, the collection of such geodesics $[ \hat{\ell}_{1}, \hat{\ell}_{2}]$ is finite by the finiteness of $L_{1}$ and
$L_{2}$, and uniqueness of geodesics in $\mathcal{L}(X)$).  This is a contradiction, because each element of $\{ \gamma_{1} \cdot \ell_{3}, \gamma_{2} \cdot \ell_{3},
\ldots, \gamma_{n} \cdot \ell_{3}, \ldots \}$ is such an $\ell$.
We conclude that each $\ell_{3} \in  [ \ell_{1},  \ell_{2}]$ has a discrete orbit.

We now show that $\ell_{3}$ is the axis of some isometry, using the fact that $\ell_{3}$ has a discrete orbit and that $\ell_{3}$ is parallel to an axis $\ell_{1}$.
We write $\ell_{1} = \{ y_{1} \} \times \mathbb{R}$ and $\ell_{3} = \{ y_{3} \} \times \mathbb{R}$.  Since $\ell_{1}$ is an axis, there is some $\gamma$ acting on 
$\ell_{1}$ by translation, and this $\gamma$ preserves the factorization of $X_{\ell_{1}} \subseteq X$, which, by definition is the union of all lines parallel to
$\ell_{1}$.  Therefore $\gamma (\hat{y}, t) = (\gamma_{1}(\hat{y}), \gamma_{2}(t))$, where $\gamma_{1}: Y \rightarrow Y$ fixes $y_{1}$, and $\gamma_{2}: \mathbb{R}
\rightarrow \mathbb{R}$ is defined by $\gamma_{2}(t) = t+K$, for some $K \in \mathbb{R}$.  Since $\ell_{3}$ has a discrete orbit and $Y$ is proper by Theorem \ref{decomp}(2), the orbit of  $y_{3}$ under the
action of $\langle \gamma_{1} \rangle$ must be finite.  Therefore, there is some power $N$ such that $\gamma_{1}^{N}(y_{3}) = y_{3}$.  This $\gamma_{1}^{N}$ acts
on $\ell_{3}$ by translation.
\end{proof}

\begin{example} \label{silly}
We now give an example of a proper CAT(0) space $X$ and an isometric action by $\Gamma$ on $X$ which demonstrates that the space $\mathcal{A}(X; \Gamma)$ can fail
to be a disjoint union of proper CAT(0) spaces.

Let $T$ be a rooted infinite binary tree.  We say that the root of the tree $T$ lies at the \emph{$0$th level} of $T$.  If $v$ is a node in the tree at the
$n$th level, then its children lie at the $(n+1)$st level.  We assign edge lengths to the tree $T$ as follows:  an edge connecting a node of level $n-1$ to one of
level $n$ has length $2^{-n}$.  The space $T$ is CAT(0) with respect to its length metric (i.e., the metric $d$ such that $d(x,y)$ is the length of the shortest path
from $x$ to $y$).  We compactify $T$ by adding all points at infinity.  The resulting space $\bar{T}$ is compact and CAT(0); each point in 
$\partial T = \bar{T} - T$ is precisely $1$ unit from the root.

Every graph automorphism $\phi: T \rightarrow T$ induces an isometry of $\bar{T}$.  We choose some automorphism $\phi: T \rightarrow T$ of infinite order.
(One way to produce such an automorphism is as follows.  Identify the set of nodes at level $n$ with the set $\mathbb{Z}/ 2^{n}\mathbb{Z}$ of left cosets of
$2^{n}\mathbb{Z}$ in $\mathbb{Z}$.  An edge connects a node $a + 2^{n+1}\mathbb{Z}$ at level $n+1$ to a node $b + 2^{n}\mathbb{Z}$ at level $n$ if and only
if $a + 2^{n+1}\mathbb{Z} \subseteq b + 2^{n}\mathbb{Z}$.  With these identifications, addition by $1$ induces an automorphism $\phi_{1}: T \rightarrow T$
of infinite order.)
We let $\hat{\phi}: \bar{T} \rightarrow \bar{T}$ denote the unique extension of $\phi: T \rightarrow T$ to the boundary.

Consider the action of $\mathbb{Z} = \langle \gamma \rangle$ on $\bar{T} \times \mathbb{R}$ defined by $\gamma \cdot (x,t) = ( \hat{\phi}(x), t+1)$.
This is a proper, cocompact, isometric action of $\mathbb{Z}$ on the proper CAT(0) space $\bar{T} \times \mathbb{R}$.  The space of lines 
$\mathcal{L}(\bar{T} \times \mathbb{R})$ is $\{ \{x \} \times \mathbb{R} \mid x \in \bar{T} \}$; this space can be identified with
$\bar{T}$.  All of the lines $\{ x \} \times \mathbb{R}$, for $x \in T$, are axes for subgroups $\langle \gamma^{N} \rangle$ (where $N$ varies with $x$).  If
$x \in \partial T$ has an infinite orbit with respect to $\langle \hat{\phi} \rangle$, then $\{ x \} \times \mathbb{R}$ isn't an axis.  It follows that
$\mathcal{A}( \bar{T} \times \mathbb{R}; \langle \gamma \rangle)$ isn't closed as a subset of $\mathcal{L}(X)$, so it cannot be proper.
\end{example}

\begin{definition} \label{wb}
Let $X$ be a proper CAT(0) space, and let $\Gamma$ act properly on $X$ by semisimple isometries.  We say that $\mathcal{A} \subseteq \mathcal{A}(X; \Gamma)$
is a \emph{well-behaved space of axes for $\Gamma$} (or $\mathcal{A}$ is a \emph{$\Gamma$-well-behaved space of axes}) 
if $\mathcal{A}$ is closed and component-wise convex as a subset of $\mathcal{L}(X)$, 
$\Gamma$-equivariant, and, whenever
a component $C$ of $\mathcal{L}(X)$ contains an element of $\mathcal{A}(X; \Gamma)$, $C$ also contains an element of $\mathcal{A}$. 
\end{definition}

\begin{remark} Note that the space $\bar{T} \times \mathbb{R}$ has many $\Gamma$-well-behaved spaces of axes.  For instance, the root of $T$ crossed with
$\mathbb{R}$ is one such well-behaved space of axes.  I don't know of a proper CAT(0) space $X$ admitting a proper isometric action by a group $\Gamma$ such
that $\Gamma$ has no well-behaved space of axes.
\end{remark}
 
\section{Good Covers for Proper CAT($0$) $\Gamma$-spaces}
 
\begin{definition} \label{goodcov}
Let $X$ be a proper CAT(0) space or disjoint union of countably many proper CAT(0) spaces, 
endowed with an isometric $\Gamma$-action having discrete orbits.  A \emph{$\Gamma$-good} cover $\mathcal{U}$ of $X$
is a cover of $X$ by open balls such that:
\begin{enumerate}
\item $\Gamma \cdot \mathcal{U} = \mathcal{U}$;
\item if $B \in \mathcal{U}$ and $\gamma \cdot B \cap B \neq \emptyset$, then $\gamma \cdot B = B$;
\end{enumerate}
\end{definition}

\begin{proposition} \label{goodex}
If $X$ is a proper CAT(0) space or a disjoint union of countably many proper CAT(0) spaces endowed with an isometric 
$\Gamma$-action having discrete orbits, then $X$ has a locally finite $\Gamma$-good cover $\mathcal{U}$.
\end{proposition}

\begin{proof}
We sketch the proof for the case in which $X$ is a proper CAT($0$) space.  If $X$ is a disjoint union of countably many proper
CAT($0$) spaces, then essentially the same argument works with only minor changes.

For each $x \in X$, we choose $0 < \epsilon_{x} < 1/2$ so that, for any $\gamma \in \Gamma$,
$\gamma \cdot B_{\epsilon_{x}}(x) \cap B_{\epsilon_{x}}(x) = \emptyset$ if $\gamma \cdot x \neq x$.  We let
$W_{x} = \bigcup_{\gamma \in \Gamma} \gamma \cdot B_{\epsilon_{x}}(x)$ and let $\mathcal{W} = \{ W_{x} \mid x \in X \}$.

Choose a basepoint $\ast \in X$.
We inductively define a subcover $\mathcal{W}_{\infty} \subseteq \mathcal{W}$ as follows.  We let $\mathcal{W}_{1} \subseteq \mathcal{W}$ 
be a finite cover of $\overline{B}_{1}(\ast)$.  (This closed ball is compact by properness of $X$.)  We let $\mathcal{W}_{2} \subseteq \mathcal{W}$
be a finite cover of $\overline{B}_{2}(\ast)$.  If $\mathcal{W}_{i}$ has been defined, then we let $\mathcal{V}_{i} = \mathcal{W}_{1} \cup \ldots \cup \mathcal{W}_{i}$.
Now suppose that $\mathcal{W}_{n} \subseteq \mathcal{W}$ has been defined, where $n \geq 2$.  We let
$\mathcal{W}_{n+1} \subseteq \mathcal{W}$ be a finite cover of 
$$\overline{B}_{n+1}(\ast) - \bigcup_{W_{x} \in \mathcal{V}_{n}} W_{x},$$
where $\mathcal{W}_{n+1}$ is chosen to be disjoint from $\overline{B}_{n-1}(\ast)$.  (Any element of $\mathcal{W}$ meeting $\overline{B}_{n-1}(\ast)$ must be completely
covered by the collection $\mathcal{V}_{n}$, so it is possible to satisfy the latter condition.)
This completes the inductive definition of $\mathcal{W}_{n}$.  We let 
$$\mathcal{W}_{\infty} = \bigcup_{n=1}^{\infty} \mathcal{W}_{n}.$$ 
Finally, we replace the cover $\mathcal{W}_{\infty}$, which consists of orbits of balls, by the individual constituent balls.  The result is locally finite,
and a $\Gamma$-good cover.
\end{proof}

\begin{remark}
One can probably control the dimension of the nerve of $N(\mathcal{U})$ (see Definition \ref{nofc} below) with sufficient care.  We don't need to do this.  In fact, 
Hatcher (\cite{Hatcher}, page 459, Corollary 4G.3) proves Lemma \ref{hatcher} without even the assumption that $\mathcal{U}$ is locally finite.
\end{remark}

\begin{definition} \label{nofc} 
Let $\mathcal{U}$ be a cover of a topological space $X$.  The \emph{nerve} $N(\mathcal{U})$ of the cover $\mathcal{U}$ is a simplicial
complex on the set $\mathcal{U}$.  A finite collection $\{ U_{1}, U_{2}, U_{3}, \ldots, U_{n} \} \subseteq \mathcal{U}$ is a simplex of
$N(\mathcal{U})$ if and only if $U_{1} \cap \ldots \cap U_{n} \neq \emptyset$.
\end{definition}

\begin{lemma} \label{hatcher} (\cite{Hatcher}, page 459, Corollary 4G.3)
If $\mathcal{U}$ is an open cover of a paracompact space $X$ such that every non-empty intersection of finitely many sets in
$\mathcal{U}$ is contractible, then $X$ is homotopy equivalent to the nerve $N(\mathcal{U})$. \qed
\end{lemma}

The first part of the following Proposition was essentially proved by Pedro Ontaneda in \cite{Pedro} (Proposition A, 2.1).  The second part
of Proposition \ref{approx} is a straightforward adaptation of his argument.

\begin{proposition} \label{approx}
Let $X$ be a proper CAT(0) space or a disjoint union of proper CAT(0) spaces.  Suppose that $\Gamma$ acts isometrically on $X$ with discrete orbits.
Let $\mathcal{U}$ be a $\Gamma$-good cover for $X$.
\begin{enumerate}
\item If $X$ is a proper CAT(0) space and $\Gamma$ acts properly on $X$, then $N(\mathcal{U})$ is an $E_{\mathcal{FIN}}\Gamma$-complex.
\item Suppose $X$ is a disjoint union of proper CAT(0) spaces.  Suppose that, for each $\gamma \in \Gamma$ of infinite order, $Fix_{X}(\gamma)$ is path connected (and thus non-empty).  Suppose that the stabilizer of each $x \in X$ is in $\mathcal{VC}$. The nerve $N(\mathcal{U})$ is an $I_{\mathcal{VC}_{\infty}}\Gamma$ complex.   
\end{enumerate}
\end{proposition}

\begin{proof}
If $\mathcal{U}$ is a $\Gamma$-good cover for $X$, 
then it is fairly clear that an element $\gamma \in \Gamma$ leaves a simplex $c \subseteq N(\mathcal{U})$ invariant if and only if $\gamma$
fixes $c$ pointwise.  We concentrate here on establishing the other conditions from Definition \ref{cx}.

\begin{enumerate}
\item Every non-empty intersection of open balls in $\mathcal{U}$ is a convex subset of $X$, and
therefore a CAT(0) space by Proposition \ref{dump}(2).  It now follows
from Lemma \ref{hatcher} that $N(\mathcal{U})$ is homotopy equivalent to $X$.

Let us suppose that $H$ is an infinite subgroup of $\Gamma$.  The fixed set $Fix_{N(\mathcal{U})}(H)$ is the full subcomplex
on the set $\{ U \mid h \cdot U = U, ~for~all~ h \in H \}$.  We claim that the latter set is empty.  For suppose $h \cdot U = U$,
for all $h \in H$.  Since $U$ has compact closure and $H$ is infinite, this violates the properness of the action of $\Gamma$.  Therefore
$Fix_{N(\mathcal{U})}(H) = \emptyset$.

Now suppose that $H$ is a finite subgroup of $\Gamma$.  The space $Fix_{X}(H)$ is a nonempty closed convex subspace of $X$ (by Corollary 2.8 on page 179 of \cite{BH}), 
and therefore
a proper CAT(0) space.  We let $\mathcal{U}_{res} = \{ U \cap Fix_{X}(H) \mid U \in \mathcal{U} ~and~ U \cap Fix_{X}(H) \neq \emptyset \}$.
The nerve $N(\mathcal{U}_{res})$ is homotopy equivalent to $Fix_{X}(H)$ for exactly the reasons that $N(\mathcal{U})$ is homotopy equivalent
to $X$.  It follows that $N(\mathcal{U}_{res})$ is contractible.

Let's consider the subcomplex $Fix_{N(\mathcal{U})}(H)$ for our finite $H$.  We claim that 
$Fix_{N(\mathcal{U})}(H)^{0} = \{ U \in \mathcal{U} \mid U \cap Fix_{X}(H) \neq \emptyset \}$ (and therefore $Fix_{N(\mathcal{U})}(H)$ is the
full subcomplex of $N(\mathcal{U})$ on the latter set).  Suppose that $U \in \mathcal{U}$ and $U \cap Fix_{X}(H) \neq \emptyset$.  Let
$x \in Fix_{X}(H) \cap U$, and let $h \in H$.  Since $h \cdot x = x$, we have that $x \in h \cdot U \cap U \neq \emptyset$, so $h \cdot U = U$.
Therefore $\{ U \in \mathcal{U} \mid U \cap  Fix_{X}(H) \neq \emptyset \} \subseteq Fix_{N(\mathcal{U})}(H)^{0}$.

Now suppose that $U \in \mathcal{U}$ and $U \cap Fix_{X}(H) = \emptyset$.  It follows that some $h \in H$ moves the center $c$ of $U$: $h \cdot c \neq c$.
This implies that $h \cdot U \neq U$, so $U \not \in Fix_{N(\mathcal{U})}(H)^{0}$.  This proves the claim.

We define a simplicial map $\rho: Fix_{N(\mathcal{U})}(H) \rightarrow N(\mathcal{U}_{res})$ on vertices, sending $U$ to $U \cap Fix_{X}(H)$.  
We claim that $\rho$ is a homotopy equivalence, because, indeed, the  inverse image of each simplex is a simplex (\cite{Pedro}, pages 50-51 uses this principle).  
We prove the latter statement.  
Let $\{ U_{1} \cap Fix_{X}(H), U_{2} \cap Fix_{X}(H), \ldots, U_{n} \cap Fix_{X}(H) \}$ be a simplex of $N(\mathcal{U}_{res})$.  Thus, in particular
$$ (\ast) \quad \bigcap_{i=1}^{n} \left( U_{i} \cap Fix_{X}(H) \right) \neq \emptyset.$$
We note that $\rho^{-1} \{ U_{1} \cap Fix_{X}(H), \ldots, U_{n} \cap Fix_{X}(H) \} = \{ U \in \mathcal{U} \mid U \cap Fix_{X}(H) = U_{i} \cap Fix_{X}(H), ~for~some~
i \in \{ 1, \ldots, n \} \}$.  (This is immediate from the definition of $\rho$.)  The latter set is indeed a simplex of $Fix_{N(\mathcal{U})}(H)$
by ($\ast$).  This proves the claim.  It now follows that $Fix_{N(\mathcal{U})}(H)$ is contractible.  This completes the proof of  (1).

\item Suppose that $H \leq \Gamma$ is not virtually cyclic.  We have that $Fix_{N(\mathcal{U})}(H)^{0} = \{ U \mid h \cdot U = U, ~for~all~ h \in H \}$.
Let us suppose that $U \in Fix_{N(\mathcal{U})}(H)^{0}$.  We let $c$ be the center of $U$.  Since $H$ is not virtually cyclic, there is some $h \in H$ such that
$h \cdot c \neq c$.  It follows that $h \cdot U \neq U$ (since $\mathcal{U}$ is a $\Gamma$-good cover, and $h \cdot U = U$ if and only if $h$ fixes $c$).  This is a contradiction, so $Fix_{N(\mathcal{U})}(H)^{0} = \emptyset$.  Therefore $Fix_{N(\mathcal{U})}(H) = \emptyset$.

Now suppose that $H \leq \Gamma$ is infinite virtually cyclic.  Let $\langle \alpha \rangle$ be a cyclic normal subgroup of finite index in $H$.  By our assumptions,
$Fix_{X}(\alpha)$ is path-connected.  In particular, $Fix_{X}(\alpha) \subseteq Y$, where $Y$ is a (proper, CAT(0)) path component of $X$.  Since 
$\langle \alpha \rangle \triangleleft H$, $H$ must leave the fixed set $Fix_{X}(\alpha)$ invariant.  It follows that $H$ acts isometries on $Y$.  Let
$h_{1}\langle \alpha \rangle$, $h_{2}\langle \alpha \rangle$, $\ldots$, $h_{m}\langle \alpha \rangle$ be a complete list of the left cosets of $\langle \alpha \rangle$
in $H$.  Let $\ast \in Fix_{X}(\alpha)$.  The orbit of $\ast$ under $H$ is $\{ h_{1} \cdot \ast, \ldots,  h_{m} \cdot \ast \}$.  The center of this orbit 
invariant under all of $H$ by Proposition \ref{dump}(3).  In particular, $Fix_{X}(H)$ is non-empty.

We summarize:
$$ Fix_{X}(H) = \bigcap_{h \in H} Fix_{X}(h)$$
is a non-empty subset of the path component $Y$, and $Fix_{X}(H)$ is therefore a closed convex subset of $Y$.  In particular, it is a proper CAT(0) space
and therefore contractible.

We let $\mathcal{U}_{res} = \{ U \cap Fix_{X}(H) \mid U \in \mathcal{U} ~and~ U \cap Fix_{X}(H) \neq \emptyset \}$ (exactly as in (1)).  As in (1), $N(\mathcal{U}_{res})$
is homotopy equivalent to $Fix_{X}(H)$, so $N(\mathcal{U}_{res})$ is contractible.

One can show that $Fix_{N(\mathcal{U})}(H)^{0} = \{ U \in \mathcal{U} \mid U \cap Fix_{X}(H) \neq \emptyset \}$ exactly as in (1).
It follows that $Fix_{N(\mathcal{U})}(H)$ is the full
subcomplex of $N(\mathcal{U})$ on  the latter set.  

We define a simplicial map $\rho: Fix_{N(\mathcal{U})}(H) \rightarrow N(\mathcal{U}_{res})$ on vertices exactly as before, sending $U$ to $U \cap Fix_{X}(H)$.
This map is a homotopy equivalence for the same reasons as in (1).  It follows that $Fix_{N(\mathcal{U})}(H)$ is contractible. 
\end{enumerate}
\end{proof}

\section{A Construction of $E_{\mathcal{VC}}(\Gamma)$ for CAT($0$) Groups}

\begin{theorem} \label{main}
Let $X$ be a proper CAT(0) space, let $\Gamma$ be a group acting properly on $X$ by semisimple isometries, and let $\mathcal{A}$ be a $\Gamma$-well-behaved
collection of axes in $X$.  There are $\Gamma$-good covers $\mathcal{U}$ of $X$, and $\mathcal{V}$ of $\mathcal{A}$.  The space
$N(\mathcal{U}) \ast N(\mathcal{V})$ is an $E_{\mathcal{VC}}\Gamma$-complex.  
\end{theorem}

\begin{proof}
By Proposition \ref{goodex}, $X$ has a $\Gamma$-good cover $\mathcal{U}$.  By Proposition \ref{goodex} 
and Definition \ref{wb}, $\mathcal{A}$ has a $\Gamma$-good cover $\mathcal{V}$.
By Proposition \ref{approx}(1), $N(\mathcal{U})$ is an $E_{\mathcal{FIN}}\Gamma$-complex.

We want to  apply Proposition \ref{approx}(2).  
Thus, we check that, for each $\gamma \in \Gamma$ of infinite order, $Fix_{\mathcal{A}}(\gamma)$ is path connected, and that the
stabilizer of each $\ell \in \mathcal{A}$ is virtually cyclic.  So suppose first that $\gamma \in \Gamma$ has infinite order.  Since $\Gamma$ acts on $X$ by semisimple
isometries, there is some line $\ell \subseteq X$ on which $\gamma$ acts by translation.  Therefore $\gamma \cdot \ell = \ell$.  Since $\mathcal{A}$ is a $\Gamma$-well-behaved space of axes, there is some $\ell_{1} \in \mathcal{A}$ parallel to $\ell$.  The orbit of $\ell_{1}$ under the  action of $\langle \gamma \rangle$
is bounded, since $d(\ell_{1}, \ell) < \infty$ and all translates of $\ell_{1}$ under $\langle \gamma \rangle$ are equidistant from $\ell$.  Since 
$\langle \gamma \rangle \cdot \ell_{1}$ is bounded, there is a center $\ell_{c} \in \mathcal{A}$ of $\langle \gamma \rangle \cdot \ell_{1}$  which is invariant under  the action
of $\langle \gamma \rangle$ (Proposition \ref{dump}(3)).  
This shows that $Fix_{\mathcal{A}}(\gamma)$ is non-empty.  Now since any two axes of $\gamma$ are parallel, $Fix_{\mathcal{A}}(\gamma)$ is
contained in a  single connected component of $\mathcal{L}(X)$, namely $\mathcal{L}_{\ell_{c}}(X)$.  We take two axes 
$\hat{\ell}_{1}, \hat{\ell}_{2} \in Fix_{\mathcal{A}}(\gamma)$.  Since $\mathcal{A}$ is a convex subset of $\mathcal{L}(X)$, the geodesic segment 
$[\hat{\ell}_{1}, \hat{\ell}_{2}]$ lies in $\mathcal{A}$.  Since  the endpoints $\hat{\ell}_{1}, \hat{\ell}_{2}$ are fixed by $\gamma$ and 
$\mathcal{L}_{\ell_{c}}(X)$ is uniquely geodesic, it must be that $[ \hat{\ell}_{1}, \hat{\ell}_{2}] \subseteq Fix_{\mathcal{A}}(\gamma)$.  This shows that 
$Fix_{\mathcal{A}}(\gamma)$ is path connected.

Now let $\ell \in \mathcal{A}$ be arbitrary.  We consider the stabilizer $\Gamma_{\ell} = \{ \gamma \in \Gamma \mid \gamma \cdot \ell = \ell \}$.  The homomorphism
$\rho: \Gamma_{\ell} \rightarrow Isom(\ell)$ must have finite kernel by the properness of the action of $\Gamma$ on $X$.  Since the orbit of a point $x \in \ell$ is
discrete, it must be that $\rho(\Gamma_{\ell})$ is isomorphic to  $\{ 1 \}$, $D_{\infty}$ (the infinite dihedral group), or $\mathbb{Z}$.  It follows that $\Gamma_{\ell}$ fits into an exact sequence
$$ F \rightarrow \Gamma_{\ell} \rightarrow C,$$
where $F$ is finite and $C$ is isomorphic to $\{ 1 \}$, $D_{\infty}$, or $\mathbb{Z}$.  
It follows that $\Gamma_{\ell}$ is virtually cyclic.

It now follows from  Proposition \ref{approx}(2) that $N(\mathcal{V})$ is an $I_{\mathcal{VC}_{\infty}}\Gamma$-complex.  Proposition \ref{promoting} implies that
$N(\mathcal{U}) \ast N(\mathcal{V})$ is an $E_{\mathcal{VC}}\Gamma$-complex.
\end{proof}

\section{A Construction of $E_{\mathcal{FBC}}(\Gamma)$ for CAT(0) Groups}

\begin{definition} \label{directed} Let $c, c': \mathbb{R} \rightarrow X$ be geodesic lines
in a CAT(0) space $X$.  We write $c \sim c'$ if there is a constant $K$
such that $c'(t) = c(t+K)$.  An equivalence class $[c]$ is called a
\emph{directed geodesic line}.

If $\mathcal{A}$ is an arbitrary subset of the space of lines 
$\mathcal{L}(X)$, we set
$$ \mathcal{DA} = \{ [c] \mid c(\mathbb{R}) \in \mathcal{A} \}.$$
We metrize $\mathcal{DA}$ as follows:
$$ d([c_{1}], [c_{2}]) = inf \{ d(c'_{1}, c'_{2}) \mid c'_{i} \sim c_{i} ~for~
i =1,2 \}.$$
It is not difficult to see that $d$ is indeed a metric on $\mathcal{DA}$.
\end{definition}

\begin{remark}  Another way to describe the metric $d$ is as follows:
The distance between the directed lines $[c_{1}], [c_{2}]$ is the same as
the distance between $c_{1}(\mathbb{R})$ and $c_{2}(\mathbb{R})$, if the latter
distance is finite (in the sense of Definition \ref{sol}) and if $c_{1}$ and $c_{2}$
are asymptotic.  Otherwise, the distance is infinite.
\end{remark}

\begin{lemma} \label{pancakes}  Let $\mathcal{A} \subseteq \mathcal{L}(X)$ be a $\Gamma$-well-behaved
space of axes.  The map $\pi: \mathcal{DA} \rightarrow \mathcal{A}$, given by sending
$[c]$ to $c(\mathbb{R})$, is a $2$-to-$1$ covering map and $\mathcal{DA}$ is a stack-of-pancakes
covering space, i.e., $\mathcal{DA}$ is isometric to the product 
$\mathcal{A} \times \{ 0, 1 \}$, where $\{ 0, 1 \}$ is given the metric in which $d(0,1) = \infty$.
In particular, $\mathcal{DA}$ is a disjoint union of proper CAT(0) spaces.

There is an involution $i: \mathcal{DA} \rightarrow \mathcal{DA}$ sending $[c]$ to $[-c]$, where
$-c(t) = c(-t)$.  The map $i$ is a covering transformation of the covering map $\pi: \mathcal{DA}
\rightarrow \mathcal{A}$.  Moreover, for any $\gamma \in \Gamma$, 
$\gamma \cdot (i([c])) = i( \gamma \cdot [c])$.

Let $\mathcal{U}$ be a $\Gamma$-good cover of $\mathcal{A}$.  If $B \in \mathcal{U}$, we let
$B^{+}$ and $B^{-}$ denote the two sheets over $B$.  The collection
$$\overline{\mathcal{U}} = \bigcup_{B \in \mathcal{U}} \{ B^{+}, B^{-} \} $$
is a $\Gamma$-good cover of $\mathcal{DA}$, and the involution $i$ acts cellularly on
$N(\overline{\mathcal{U}})$.
\end{lemma}

\begin{proof}
The group $\langle i \rangle \cong \mathbb{Z}/2\mathbb{Z}$ acts freely and isometrically on
$\mathcal{DA}$.  Since the group $\langle i \rangle$ is finite, the action is therefore by 
covering transformations.  The quotient of $\mathcal{DA}$ by the action of $\langle i \rangle$
is $\mathcal{A}$.  We let $\pi: \mathcal{DA} \rightarrow \mathcal{A}$ denote the quotient 
projection.  The map $\pi$ is a $2$-to-$1$ covering map.  Since each connected component of
$\mathcal{A}$ is simply connected (indeed, contractible), $\mathcal{DA}$ is a stack-of-pancakes
covering space.

It is straightforward to check that $\gamma \cdot (i ([c])) = i ( \gamma \cdot [c])$, and that
$\overline{\mathcal{U}}$ is a $\Gamma$-good cover of $\mathcal{DA}$.  
\end{proof}

In what follows, we equip the infinite-dimensional sphere $\mathbb{S}^{\infty}$ with a cell structure such that the antipodal map 
$a: \mathbb{S}^{\infty} \rightarrow \mathbb{S}^{\infty}$ acts in a cell-permuting way, and $a \cdot c \cap c = \emptyset$ for each cell
$c \subseteq \mathbb{S}^{\infty}$.

\begin{theorem} \label{biggie}
Suppose that $\mathcal{A}$ is a $\Gamma$-well-behaved collection of axes, $\mathcal{U}$ is a
$\Gamma$-good cover of $\mathcal{A}$, and let $\overline{\mathcal{U}}$ be the $\Gamma$-good cover
of $\mathcal{DA}$ determined by the previous lemma.  

The space 
$$K = \left( N(\overline{\mathcal{U}}) \times \mathbb{S}^{\infty} \right)/ \sim.$$
is an $I_{\mathcal{FBC}_{\infty}}\Gamma$-complex, where $(x,p) \sim (i(x), a(p))$, and, for
any $\gamma \in \Gamma$, $\gamma \cdot [(x,p)] = [(\gamma \cdot x, p)]$.
\end{theorem}

\begin{proof}
We first note that $K$ has a natural cell structure, since it is the quotient of the CW-complex
$N(\overline{\mathcal{U}}) \times \mathbb{S}^{\infty}$ by a free cell-permuting action of
$\mathbb{Z}/ 2\mathbb{Z} \cong \langle t \rangle$, where $t \cdot (x,p) = (i(x), a(p))$.

Now we show that if a cell of $K$ is left invariant as a set by some $\gamma \in \Gamma$, 
then the cell is actually fixed pointwise by $\gamma$.  A cell $e$ of $K$ has the form
$$ \left[ (c_{1} \times c_{2}) \cup ( i(c_{1}) \times a(c_{2}) ) \right] / \sim.$$
If this cell is invariant under the action of $\gamma$, it must be that the cell
$\pi(c_{1}) \in N(\mathcal{U})$ (where $\pi$ is induced by the projection 
$\pi: \mathcal{DA} \rightarrow \mathcal{A}$) is fixed pointwise by the action of $\gamma$.
It follows that the action of $\gamma$ on $c_{1}$ either agrees with that of $i$, or 
$\gamma_{\mid c_{1}} = id_{\mid c_{1}}$.  We choose $(x,y) \in c_{1} \times c_{2}$ to represent a
typical element of $e$.  If $\gamma$ agrees with $i$ on $c_{1}$, then 
$\gamma \cdot (x,y) = (i(x),y)$.  This point is not an element of $e$.  This is a contradiction,
and shows that $\gamma$ agrees with the identity on $c_{1}$.  Therefore $\gamma \cdot (x,y) = (x,y)$,
as we wished to show.

The same line of reasoning establishes the identity
$$Fix_{K}(H) = 
\left( Fix_{N(\overline{\mathcal{U}})}(H) \times \mathbb{S}^{\infty} \right)/ \sim.$$

We verify properties (1) and (2) of $I_{\mathcal{FBC}_{\infty}}\Gamma$-complexes from Definition \ref{cx}.
Let $H \in \mathcal{FBC}_{\infty}$.  We use the following fact about infinite
finite-by-cyclic groups $H$:  if $H$ acts properly and isometrically on a line $\ell$, then the
image of the associated homomorphism $\rho: H \rightarrow Isom(\ell)$ is a group of translations.
Since $H$ is infinite and virtually cyclic, it follows from Proposition 3.9(2) that
$Fix_{N(\mathcal{U})}(H)$ is a contractible subcomplex of $N(\mathcal{U})$.

We claim that $\pi^{-1}(Fix_{N(\mathcal{U})}(H)) = Fix_{N(\overline{\mathcal{U}})}(H)$.  Suppose that
$c$ is a cell of $N(\overline{\mathcal{U}})$, and $c \subseteq \pi^{-1}(Fix_{N(\mathcal{U})}(H))$.
We have that $\pi(c)$ is fixed by all of $H$, so each $h \in H$ either acts as the identity on $c$,
or as the inversion $i$.  However, if $h \in H$ acts as inversion on $c$, then $h$ must act as the
inversion on some line $\ell \subseteq X$.  This contradicts the fact that 
$\rho(H) \subseteq Isom(\ell)$ must be a group of translations.  It follows that each $h \in H$
fixes $c$, so $\pi^{-1}(Fix_{N(\mathcal{U})}(H)) \subseteq Fix_{N(\overline{\mathcal{U}})}(H)$.  The
reverse inclusion is clear.

It now follows that 
$$ Fix_{K}(H) = \left( \pi^{-1}(Fix_{N(\mathcal{U})}(H)) \times \mathbb{S}^{\infty} \right)/\sim.$$
Thus, $Fix_{K}(H)$ is the quotient of two contractible sets by an involution $(x,y) \rightarrow
(i(x), a(y))$ which identifies these components by a homeomorphism.  Therefore, $Fix_{K}(H)$
is contractible.  

If $H \in \mathcal{VC} - \mathcal{FBC}$, then whenever $H$ acts properly and isometrically on a line
$\ell$, the image of the associated homomorphism $\rho: H \rightarrow Isom(\ell)$ has an involution.
The group $H$ has a subgroup of index two, $H^{+}$, where $H^{+} \in \mathcal{FBC}$.  We let $\alpha$
be an element of order $2$ such that $H = \langle H^{+}, \alpha \rangle$.  By the above argument,
$\pi^{-1}(Fix_{N(\mathcal{U})}(H^{+}) = Fix_{N(\overline{\mathcal{U}})}(H^{+})$, and this space is a
disjoint union of two contractible pieces which are swapped by the covering transformation $i$.
We note that the element $\alpha$ must act as an inversion on any line $\ell$ on which $H$ acts
properly and isometrically, for otherwise it would necessarily act as the identity, and then the
entire group $H$ would map to a group of translations of $\ell$ under $p: H \rightarrow Isom(\ell)$,
and this is impossible.  It follows that $\alpha$ swaps the components of 
$\pi^{-1}(Fix_{N(\mathcal{U})}(H^{+}))$; in particular
$Fix_{N(\overline{\mathcal{U}})}(H) = \emptyset$, so $Fix_{K}(H) = \emptyset$.

We note finally that $Fix_{N(\overline{\mathcal{U}})}(H) = \emptyset$ if $H \not \in \mathcal{VC}$,
and so $Fix_{K}(H) = \emptyset$ if $H \not \in \mathcal{VC}$.  This completes the proof.
\end{proof}

\begin{corollary} \label{FBC} If $\Gamma$ acts properly and isometrically by semisimple isometries on the
proper CAT(0) space $X$, $\mathcal{A}$ is a $\Gamma$-well-behaved space of axes, $\mathcal{U}$
is a $\Gamma$-good cover of $X$, and $\mathcal{V}$ is a $\Gamma$-good cover of $\mathcal{A}$, then
the space
$$ N(\mathcal{U}) \ast K $$
is a $E_{\mathcal{FBC}}\Gamma$-complex, where $K$ is as defined in the previous theorem.
\end{corollary}

\begin{proof}
This follows immediately from the Theorem \ref{biggie} and from Proposition \ref{promoting}.
\end{proof}

\bibliography{EVCfinrev}

\begin{thebibliography}{1}

\bibitem{Bartels}
Arthur~C. Bartels.
\newblock On the domain of the assembly map in algebraic {$K$}-theory.
\newblock {\em Algebr. Geom. Topol.}, 3:1037--1050 (electronic), 2003.

\bibitem{BH}
Martin~R. Bridson and Andr{\'e} Haefliger.
\newblock {\em Metric spaces of non-positive curvature}, volume 319 of {\em
  Grundlehren der Mathematischen Wissenschaften [Fundamental Principles of
  Mathematical Sciences]}.
\newblock Springer-Verlag, Berlin, 1999.

\bibitem{FJ}
F.~T. Farrell and L.~E. Jones.
\newblock Isomorphism conjectures in algebraic {$K$}-theory.
\newblock {\em J. Amer. Math. Soc.}, 6(2):249--297, 1993.

\bibitem{CFH}
Benjamin~Fehrman Frank~Connolly and Michael Hartglass.
\newblock On the dimension of the virtually cyclic classifying space of a
  crystallographic group.
\newblock {\em Preprint, 2006; arxiv:math.AT/0610387}.

\bibitem{Hatcher}
Allen Hatcher.
\newblock {\em Algebraic topology}.
\newblock Cambridge University Press, Cambridge, 2002.

\bibitem{Khan}
Qayum~Khan James F.~Davis and Andrew Ranicki.
\newblock Algebraic k-theory over the infinite dihedral group.
\newblock {\em Preprint, 2008; arxiv:0803.1639}.

\bibitem{Lueck}
Wolfgang Lueck.
\newblock On the classifying space of the family of finite and of virtually
  cyclic subgroups for cat(0) groups.
\newblock {\em Preprint, 2009; arxiv:0902.0718}.

\bibitem{Munk}
James~R. Munkres.
\newblock {\em Topology: a first course}.
\newblock Prentice-Hall Inc., Englewood Cliffs, N.J., 1975.

\bibitem{Pedro}
Pedro Ontaneda.
\newblock Cocompact {CAT}(0) spaces are almost geodesically complete.
\newblock {\em Topology}, 44(1):47--62, 2005.

\end{thebibliography}
\nocite{*}
\bibliographystyle{plain}

\end{document}